\documentclass{amsart}
\usepackage{amsfonts}

\newtheorem{theorem}{Theorem}
\theoremstyle{plain}

\newtheorem{corollary}{Corollary}

\newtheorem{example}{Example}

\newtheorem{lemma}{Lemma}

\newtheorem{proposition}{Proposition}

\numberwithin{equation}{section}
\input{tcilatex}

\begin{document}
\title[Abelian state-closed groups]{Abelian state-closed subgroups of
automorphisms of $m$-ary trees}
\author{Andrew M. Brunner}
\address{Department of Mathematics, University of Wisconsin-Parkside,
Kenosha, Wisconsin 53141-2000}
\email{brunner@uwp.edu}
\author{Said N. Sidki}
\address{Departamento de Matematica, Universidade de Bras\'{\i}lia\\
Bras\'{\i}lia DF 70910-900, Brazil}
\email{sidki@mat.unb.br}
\thanks{The first author is grateful for the hospitality of the Department
of Mathematics of the University of Bras\'{\i}lia as well as travel
assistance from the University of Wisconsin-Parkside. The second author
acknowledges support from the Brazilian scientific agencies CNPq and FAPDF.
Both authors are thankful to Laurent Bartholdi for making numerous comments.}
\date{October 02, 2009}
\subjclass[2000]{ Primary 20E08, 20F18}
\keywords{Automorphisms of trees, state-closed groups, self-similar groups,
abelian groups, topological closure, $p$-adic integers, pro-$p$ groups.}

\begin{abstract}
The group $\mathcal{A}_{m}$ of automophisms of a one-rooted $m$-ary tree
admits a diagonal monomorphism which we denote by $x$. Let $A$ be an abelian
state-closed (or self-similar) subgroup of $\mathcal{A}_{m}$. We prove that
the combined diagonal and tree-topological closure $A^{\ast }$ of $A$ is
additively a finitely presented $\mathbb{Z}_{m}\left[ \left[ x\right] \right]
$-module where $\mathbb{Z}_{m}$ is the ring of $m$-adic integers. Moreover,
if $A^{\ast }$ is torsion-free then it is a finitely generated pro-$m$
group. Furthermore, the group $A$ splits over its torsion subgroup. We study
in detail the case where $A^{\ast }$ is additivley a cyclic $\mathbb{Z}_{m}%
\left[ \left[ x\right] \right] $-module and we show that when $m$ is a prime
number then $A^{\ast }$ is conjugate by a tree automorphism to one of two
specific types of groups.
\end{abstract}

\maketitle

\section{Introduction}

Automorphisms of one-rooted regular trees $\mathcal{T}\left( Y\right) $
indexed by finite sequences from a finite set $Y$ of size $m\geq 2$, have a
natural interpretation as automata on the alphabet $Y$, and with states
which are again automorphisms of the tree. A subgroup of the group of
automorphisms $\mathcal{A}\left( Y\right) $ of the tree is said to be 
\textit{state-closed,} in the language of automata (or \textit{self-similar}
in the language of dynamics) of degree $m$ provided the states of its
elements are themselves elements of the same group. If the group is not
state-closed then we may consider its state-closure. The prime example of a
state-closed group is the group generated by the binary adding machine $\tau
=\left( e,\tau \right) \sigma $ where $\sigma $ is the transposition $(0,1)$.

We study in this paper representations of general abelian groups as
state-closed groups of degree $m$. For this purpose we use topological and
diagonal closure operations in the automorphism group of the tree.
Representations of free abelian groups of finite rank as state-closed groups
of degree $2$ were characterized in \cite{Nek-Sid}.

An automorphism group $G$ of the tree group is said to be \textit{transitive}
provided the permutation group $P\left( G\right) $ induced by $G$ on the set 
$Y$ is transitive; actions of groups on sets, will be applied on the right.
It will be shown that the structure of state-closed groups can in a certain
sense be reduced to those which are transitive.

The automorphism group $\mathcal{A}\left( Y\right) $ of the tree is a
topological group with respect to the topology inherited from the tree. This
topology allows us to exponentiate elements of $\mathcal{A}\left( Y\right) $
by $m$-ary integers from $\mathbb{Z}_{m}$. Given a subgroup $G$ of $\mathcal{%
A}\left( Y\right) $, its topological closure $\overline{G}$ with respect to
the tree topology belongs to the same variety as $G$. Also, if $G$ is
state-closed then so is $\overline{G}$.

The diagonal map $\alpha \rightarrow \alpha ^{\left( 1\right) }=\left(
\alpha ,\alpha ,...,\alpha \right) $ is a monomorphism of $\mathcal{A}_{m}$.
Define inductively $\alpha ^{\left( 0\right) }=\alpha ,\alpha ^{\left(
i+1\right) }=\left( \alpha ^{\left( i\right) }\right) ^{\left( 1\right) }$
for $i\geq 0$. It is convenient to introduce a symbol $x$ and write $\alpha
^{\left( i\right) }$ as $\alpha ^{x^{i}}$ for $i\geq 0$. This will permit
more general exponentiation, by formal power series $p\left( x\right) \in 
\mathbb{Z}_{m}\left[ \left[ x\right] \right] $. Given a subgroup $G$ of $%
\mathcal{A}\left( Y\right) $, its \textit{diagonal closure} is the group $%
\widetilde{G}=\left\langle G^{\left( i\right) }\mid i\geq 0\right\rangle $.
Observe that the diagonal closure operation preserves the state-closed
property.

We will show that given an abelian transitive state-closed group $A$, its
diagonal closure $\widetilde{A}$ is again abelian. The composition of the
diagonal and topological closures when applied to $A$ produces an abelian
group denoted by $A^{\ast }$ which can be viewed additively as a finitely
generated $\mathbb{Z}_{m}\left[ \left[ x\right] \right] $-module. This
approach was first used in \cite{Bru-Sid97}.

The prime decomposition $m$ $=\dprod\limits_{1\leq i\leq s}p_{i}^{k_{i}}$,
provides us with the decomposition $\mathbb{Z}_{m}=\oplus _{1\leq i\leq
s}\varepsilon _{i}\mathbb{Z}_{p_{i}^{k_{i}}}$ where $\varepsilon _{i}$ are
orthogonal idempotents such that $1=\sum_{1\leq i\leq s}\varepsilon _{i}$,
and provides us also with the decomposition $\mathbb{Z}_{m}\left[ \left[ x%
\right] \right] =\oplus _{1\leq i\leq s}\varepsilon _{i}\mathbb{Z}%
_{p_{i}^{k_{i}}}\left[ \left[ x\right] \right] $. When $m=p^{k}$ and $p$ a
prime number, the rings $\mathbb{Z}_{m}\left[ \left[ x\right] \right] $ and $%
\mathbb{Z}_{p}\left[ \left[ x\right] \right] $ are isomorphic, yet when $k>1$%
, they are different representations of the same object and for this reason
we distinguish between them.

In Sections 3 and 4 we prove

\begin{theorem}
Let $A$ be an abelian transitive state-closed group of degree $m$. Then, (1)
the group $A^{\ast }$ is isomorphic to a finitely presented $\mathbb{Z}_{m}%
\left[ \left[ x\right] \right] $-module; (2) if $A^{\ast }$is torsion-free
then it is a finitely generated $\mathbb{Z}_{m}$-module which is also a pro-$%
m$ group.
\end{theorem}

Item (1) is part of Theorem 5 and item (2) is Corollary 1 of Theorem 6.

We consider in Section 5 torsion subgroups of state-closed abelian groups
and use methods from virtual endomorphisms of groups (see \cite{Nek}, \cite%
{Ber-Sid}; reviewed in Subsection 5.1) to prove the following structural
result.

\begin{theorem}
Let $A$ be an abelian transitive state-closed group of degree $m$ and $%
tor\left( A\right) $ its torsion subgroup. Then, (i) $tor\left( A\right) $
is a direct summand of $A$ and has exponent a divisor of the exponent of $%
P\left( A\right) $; (ii) the action of $A$ on the $m$-ary tree induces
transitive state-closed representations of $tor\left( A\right) $ on the $%
m_{1}$-tree and of $\frac{A}{tor\left( A\right) }$ on the $m_{2}$-tree,
where $m_{1}=|P\left( tor\left( A\right) \right) |$ and $m_{2}=|\frac{%
P\left( A\right) }{P\left( tor\left( A\right) \right) }|$; (iii) if $%
A=tor\left( A\right) $ and $P\left( A\right) \cong \oplus _{1\leq i\leq k}%
\frac{\mathbb{Z}}{m_{i}\mathbb{Z}}$, then $A^{\ast }\cong $ $\oplus _{1\leq
i\leq k}\frac{\mathbb{Z}}{m_{i}\mathbb{Z}}\left[ \left[ x\right] \right] $.
\end{theorem}

The above results are analogous to Theorem 4.3.4 of \cite{Rib-Zal} on the
structure of finitely generated pro-$p$ groups. By item (i) of the theorem,
an abelian torsion group $G$ of infinite exponent cannot have a faithful
representation as a transitive state-closed group for any finite degree. Put
differently, the group $G$ does not admit any simple virtual endomorphism.
On the other hand, the group of automorphisms of the $p$-adic tree is
replete with abelian $p$-subgroups of infinite exponent. Item (iii) follows
from Theorem 7 which is a conjugacy result and therefore more general than
isomorphism.

We focus our attention in Section 6 on transitive state-closed abelian
groups $A$ for which $A^{\ast }$ is additively a cyclic $\mathbb{Z}_{m}\left[
\left[ x\right] \right] $-module. We show

\begin{theorem}
(1) Let $q_{\mathbf{1}},...,q_{\mathbf{m}}\in \mathbb{Z}_{m}\left[ \left[ x%
\right] \right] $ and let $\sigma $ be the cycle $\left( 1,2,...,m\right) $.
Then the expression 
\begin{equation*}
\alpha =\left( \alpha ^{q_{\mathbf{1}}},...,\alpha ^{q_{\mathbf{m}}}\right)
\sigma
\end{equation*}%
is a well-defined automorphism of the $m$-ary tree and the state-closure $A$
of $\left\langle \alpha \right\rangle $ is an abelian transitive group. The
group $A^{\ast }$ is additively isomorphic to the quotient ring $\frac{%
\mathbb{Z}_{m}\left[ \left[ x\right] \right] }{\left( r\right) }\ $where%
\begin{equation*}
r=m-xq\text{ and }q=q_{\mathbf{1}}+...+q_{\mathbf{m}}\text{.}
\end{equation*}%
(2) Let $A$ be a transitive state-closed abelian group of degree $m$ such
that $A^{\ast }$ is additively a cyclic $\mathbb{Z}_{m}\left[ \left[ x\right]
\right] $-module. Then $P\left( A\right) $ is cyclic, say generated by $%
\sigma $, and $A^{\ast }$ is the state-diagonal-topological closure of an
element of the form $\alpha =\left( \alpha ^{q_{\mathbf{1}}},...,\alpha ^{q_{%
\mathbf{m}}}\right) \sigma $ for some $q_{\mathbf{1}},...,q_{\mathbf{m}}\in 
\mathbb{Z}_{m}\left[ \left[ x\right] \right] $.
\end{theorem}

Finally, we provide a complete description of the group $A^{\ast }$ for
state-closed groups of prime degree. Let $j\geq 1$ and let $\mathit{D}%
_{m}\left( j\right) $ be the group generated by the set of states of the
generalized adding machine $\alpha =\left( e,...,e,\alpha ^{x^{j-1}}\right)
\sigma $ acting on the $m$-ary tree with $\sigma =\left( 1,2,...,m\right) $.
The topological closure of $\mathit{D}_{m}\left( j\right) $ seen as $\mathbb{%
Z}_{m}$-module, is isomorphic to the ring $\frac{\mathbb{Z}_{m}\left[ \left[
x\right] \right] }{\left( r\right) }$, $r=m-x^{j}$.

\begin{theorem}
Let $A$ be an abelian transitive state-closed group of prime degree $m$ and
let $\sigma $ be the $m$-cycle automorphism. If $tor\left( A\right) $ is
nontrivial then $A^{\ast }$ is a torsion group conjugate to $\left\langle
\sigma \right\rangle ^{\ast }\left( \cong \frac{\mathbb{Z}}{m\mathbb{Z}}%
\left[ \left[ x\right] \right] \right) $. If $A$ is torsion-free then $%
A^{\ast }$ is a torsion-free group conjugate to the topological closure of $%
\mathit{D}_{m}\left( j\right) $ for some $j$.
\end{theorem}

One of the questions which has remained unanswered is whether a free abelian
group of infinite rank admits a faithful transitive state-closed
representation, even of prime degree.

\section{Preliminaries}

We fix the notation $Y=\left\{ 1,2,...,m\right\} $, $\mathcal{T}_{m}=%
\mathcal{T}\left( Y\right) $, $\mathcal{A}_{m}=\mathcal{A}\left( Y\right) $
and we let $Perm\left( Y\right) $ be the group of permutations of $Y$. A
permutation $\gamma \in Perm\left( Y\right) $ is extended to an automorphism
of the tree by $\gamma :yu\rightarrow y^{\gamma }u$, fixing the non-initial
letters of every sequence. An automorphism $\alpha \in \mathcal{A}_{m}$ is
represented as $\alpha =\left( \alpha _{1},\alpha _{2},...,\alpha
_{m}\right) \sigma \left( \alpha \right) $ where $\alpha _{i}\in \mathcal{A}%
_{m}$ and $\sigma \left( \alpha \right) \in Perm\left( Y\right) $.
Successive developments of $\alpha _{i}$ produce for us $\alpha _{u}$ (a
state of $\alpha $) for every finite string $u$ over $Y$.

The product of $\alpha =\left( \alpha _{1},\alpha _{2},...,\alpha
_{m}\right) \sigma \left( \alpha \right) $ and $\beta =\left( \beta
_{1},\beta _{2},...,\beta _{m}\right) \sigma \left( \beta \right) $ in $%
\mathcal{A}_{m}$, is 
\begin{equation*}
\alpha \beta =\left( \alpha _{1}\beta _{\left( 1\right) \sigma \left( \alpha
\right) },...,\alpha _{m}\beta _{\left( m\right) \sigma \left( \alpha
\right) }\right) \sigma \left( \alpha \right) \sigma \left( \beta \right) 
\text{.}
\end{equation*}

Let $G$ be a subgroup of $\mathcal{A}_{m}$. Denote the subgroup of $G$ which
fixes the vertices of the $i$-th level of the tree by $Stab_{G}\left(
i\right) $. Given $y\in Y$, denote by $Fix_{G}\left( y\right) $ the subgroup
of $G$ consisting of the elements of $G$, which fix $y$. The group $G$ is
said to be \textit{recurrent} provided it is transitive and $Fix_{G}\left(
1\right) $ projects in the $1$st coordinate onto $G$.

The group $\mathcal{A}_{m}$ is the inverse limit of its quotients by the $i$%
-th level stabilizers $Stab_{\mathcal{A}_{m}}(i)$ of the tree and is as such
a topological group where each $Stab_{\mathcal{A}_{m}}(i)$ is an open and
closed subgroup. For a subgroup $G$ of automorphisms of the tree, its
topological closure $\overline{G}$ coincides with the set of all infinite
products $...g_{\mathbf{i}}...g_{\mathbf{1}}g_{\mathbf{0}}$, or alternately, 
$g_{\mathbf{0}}g_{\mathbf{1}}...g_{\mathbf{i}}...$ where $g_{\mathbf{i}}\in
Stab_{G}\left( \mathbf{i}\right) $. The group $\overline{G}$ satisfies the
same group identities as $G$. We note that the property of being
state-closed is also preserved by the topological closure operation.

Let $\alpha $ be an automorphism of the tree. Then $\overline{\left\langle
\alpha \right\rangle }=\left\{ \alpha ^{p}\mid p\in \mathbb{Z}_{m}\right\} $%
. More generally, for $q=\sum_{i\geq 0}q_{i}x^{i}\in \mathbb{Z}_{m}\left[ %
\left[ x\right] \right] $ with $q_{i}\in \mathbb{Z}_{m}$, we write the
expression 
\begin{equation*}
\alpha ^{q}=\alpha ^{q_{0}}\alpha ^{q_{1}x}...\alpha ^{q_{i}x^{i}}...
\end{equation*}%
which can be verified to be a well-defined automorphism of the tree.

We recall the reduction of group actions to transitive ones, with a view to
a similar reduction for state-closed groups of automorphisms of trees. Let $%
G $ be a subgroup of $Perm\left( Y\right) $, let $\left\{ Y_{i}\mid
i=1,...,s\right\} $ be the set of orbits of $G$ on $Y$ and let $\left\{ \rho
_{i}:G\rightarrow Perm\left( Y_{i}\right) \mid i=1,...,s\right\} $ be the
set of induced representations. Then, each $\rho _{i}$ is transitive and $%
\rho :G\rightarrow \dprod\limits_{1\leq i\leq s}Perm\left( Y_{i}\right) \leq
Perm\left( Y\right) $ defined by $g\rightarrow \left( g^{\rho
_{1}},...,g^{\rho _{s}}\right) $ is a monomorphism. The reduction for tree
actions follows from

\begin{lemma}
Let $G$ be a state-closed group of automorphisms of the tree $\mathcal{T}%
\left( Y\right) $ and let $X$ be a $P\left( G\right) $-invariant subset of $%
Y $. Then, $\mathcal{T}\left( X\right) $ is $G$-invariant and for the
resulting representation $\mu :G\rightarrow \mathcal{A}\left( X\right) $,
the group $G^{\mu }$ is state-closed. If $G$ is diagonally closed or is
topologically closed then so is $G^{\mu }$.
\end{lemma}

\begin{proof}
Let $xu$ be a sequence from $X$ and let $\alpha \in G$. Then, $\left(
xu\right) ^{\alpha }=x^{\sigma \left( \alpha \right) }u^{\alpha _{x}}$. As $%
x^{\sigma \left( \alpha \right) }\in X$ and $\alpha _{x}\in G$, it follows
that $\left( xu\right) ^{\alpha }$ is a sequence from $X$. Also, for any
sequence $u$ from $X$, we have $\left( \alpha ^{\mu }\right) _{u}=\left(
\alpha _{u}\right) ^{\mu }$. Thus, $G^{\mu }$ is state-closed. The last
assertion is clear.
\end{proof}

We note the following important properties of transitive state-closed
abelian groups $A$.

\begin{proposition}
Let $A$ be an abelian transitive state-closed group of degree $m$. Then $%
Stab_{A}\left( i\right) \leq A^{\left( i\right) }$ for all $i\geq 0$. The
group $\widetilde{A}$ is an abelian transitive state-closed group and is a
minimal recurrent group containing $A$. Moreover, the topological closure
and diagonal closure operations commute when applied to $A$. The
diagonal-topological closure $A^{\ast }$ of $A$ is an abelian transitive
state-closed group.
\end{proposition}

\begin{proof}
Let $\alpha =\left( \alpha _{1},...,\alpha _{m}\right) \sigma ,\beta =\left(
\beta _{1},...,\beta _{m}\right) \in A$. Then, the conjugate of $\beta $ by $%
\alpha $ is 
\begin{equation*}
\beta ^{\alpha }=\left( \beta _{1}^{\alpha _{1}},...,\beta _{m}^{\alpha
_{m}}\right) ^{\sigma }\text{.}
\end{equation*}%
As $\alpha _{i},\beta _{i}\in A$ and $A$ is abelian, it follows that $\beta
=\left( \beta _{1},...,\beta _{m}\right) ^{\sigma }$. Furthermore, since $A$
is transitive, $\beta =\left( \beta _{1},...,\beta _{1}\right) =\left( \beta
_{1}\right) ^{\left( 1\right) }$. Thus, $Stab_{A}\left( i\right) \leq
A^{\left( i\right) }$ for all $i$. A similar verification shows that $%
\widetilde{A}=\left\langle A^{\left( i\right) }\mid i\geq 0\right\rangle $
is abelian.

Let $G$ be a recurrent group such that $A\leq G\leq \widetilde{A}$. Given $%
\alpha \in G$, as $G$ is recurrent, there exists $\beta \in Stab_{G}\left(
1\right) $ such that $\beta =\left( \beta _{1},...,\beta _{m}\right) $ with $%
\beta _{1}=\alpha $. Since $G$ is transitive and abelian, we have $\beta
_{1}=...=\beta _{m}=\alpha $; that is, $\beta =\alpha ^{\left( 1\right) }$.
Hence, $A^{\left( i\right) }\leq G$ and $G=\widetilde{A}$ follows.

The last two assertions of the proposition are clear.
\end{proof}

The following result indicates the smallness of recurrent transitive abelian
groups, from the point of view of centralizers.

\begin{proposition}
(Theorem 7 \cite{Ber-Sid}) (1) Let $A$ be a recurrent abelian group of
degree $m$ and let $C_{\mathcal{A}_{m}}\left( A\right) $ be the centralizer
of $A$ in $\mathcal{A}_{m}$. Then, $C_{\mathcal{A}_{m}}\left( A\right) =%
\overline{A}$. (2) Let $m$ be a prime number and $A$ be an infinite
transitive state-closed abelian group. Then, $C_{\mathcal{A}_{m}}\left(
A\right) =\overline{A}$.
\end{proposition}

This result will be used in the proofs of Lemma 3 and Step 4 of Theorem 9.

\section{A presentation for $A^{\ast }$}

Let $A$ be a transitive abelian state-closed group of degree $m$ and let $%
A^{\ast }$ be its diagonal-topological closure. Then $A^{\ast }$ is
additively a $\mathbb{Z}_{m}\left[ \left[ x\right] \right] $-module having
the following properties. Given $\alpha \in A^{\ast }$, then 
\begin{equation*}
\text{(i) }x\alpha =0\text{ implies }\alpha =0;\text{(ii) }m\alpha =x\gamma 
\text{ for some }\gamma \in A^{\ast }\text{.}
\end{equation*}

Let $P\left( A\right) $ be given by its presentation 
\begin{equation*}
\left\langle \sigma _{i}\text{ }\left( 1\leq i\leq k\right) \mid \sigma
_{i}^{m_{i}}=e\text{, abelian}\right\rangle \text{.}
\end{equation*}

Choose for each $\sigma _{i}$ an element $\beta _{i}$\ in $A$, which induces 
$\sigma _{i}$ on $Y$; denote $\beta _{i}$ by $\beta \left( \sigma
_{i}\right) $.\ Then for any $n\geq 0$, the automorphism of the tree $\beta
\left( \sigma _{i}\right) ^{\left( n\right) }$ is an element of $\widetilde{A%
}$ which induces $\left( \sigma _{i}\right) ^{\left( n\right) }$ on the $%
\left( n+1\right) $-th level of the tree. Although the notation $\beta _{i}$
has been used to indicate the $i$th entry in an automorphism $\beta $, we
hope this new usage will not cause confusion.

\begin{theorem}
Let $A$ be a transitive abelian state-closed group of degree $m$. Then $%
A^{\ast }$ is additively a $\mathbb{Z}_{m}\left[ \left[ x\right] \right] $%
-module generated by%
\begin{equation*}
\left\{ \beta _{i}\mid 1\leq i\leq k\right\}
\end{equation*}%
subject to the set of defining relations 
\begin{equation*}
\left\{ r_{i}=\sum_{1\leq j\leq k}m_{i}\beta _{i}-p_{ij}\beta _{j}x=0\mid
1\leq i\leq k\right\} \text{ for some }p_{ij}\in \mathbb{Z}_{m}\left[ \left[
x\right] \right] \text{.}
\end{equation*}%
Moreover, there exist $r,q\in \mathbb{Z}_{m}\left[ \left[ x\right] \right] $
such that $r=m-xq$ and $rA^{\ast }=\left( 0\right) $. The elements of $%
A^{\ast }$ can be represented additively as $\sum_{1\leq i\leq k}p_{i}\beta
_{i}$ where $p_{i}=\sum_{j\geq 0}p_{ij}x^{j}$ and each $p_{ij}\in \mathbb{Z}$
with $0\leq p_{ij}<m$.
\end{theorem}

\begin{proof}
Let $\alpha \in A^{\ast }$ and $\sigma \left( \alpha \right)
=\dprod\limits_{1\leq i\leq k}\sigma _{i}^{r_{i1}}$, $0\leq r_{i1}<m_{i}$.
Then, either $\alpha \left( \dprod\limits_{1\leq i\leq k}\beta
_{i}^{r_{i1}}\right) ^{-1}$ is the identity element or there exists $%
l_{2}\geq 1$ such that 
\begin{equation*}
\alpha \left( \dprod\limits_{1\leq i\leq k}\beta _{i}^{r_{i1}}\right)
^{-1}\in Stab\left( l_{2}\right) \backslash Stab\left( l_{2}+1\right)
\end{equation*}%
and so, $\alpha \left( \dprod\limits_{1\leq i\leq k}\beta
_{i}^{r_{i1}}\right) ^{-1}=\left( \gamma \right) ^{\left( l_{2}\right) }$
for some $\gamma \in A^{\ast }$. We treat $\gamma $ in the same manner as $%
\alpha $. In the limit, we obtain%
\begin{eqnarray*}
\alpha &=&\dprod\limits_{1\leq i\leq k}\left( \beta _{i}^{r_{i1}}\left(
\beta _{i}^{r_{i2}}\right) ^{\left( l_{2}\right) }...\left( \beta
_{i}^{r_{ij}}\right) ^{\left( l_{j}\right) }...\right) \\
&=&\dprod\limits_{1\leq i\leq k}\beta _{i}^{q_{i}}
\end{eqnarray*}%
where $0\leq r_{ij}<m_{i}$, $1\leq l_{2}<l_{3}<..<l_{j}<...$ and where $%
q_{i}=r_{i1}+\dsum\limits_{j\geq 2}r_{ij}x^{l_{j}}$ are formal power series
in $x$. Additively then,%
\begin{equation*}
\alpha =\sum_{1\leq i\leq k}q_{i}\beta _{i}\in \sum_{1\leq i\leq k}\mathbb{Z}%
_{m_{i}}\left[ \left[ x\right] \right] \beta _{i}\text{.}
\end{equation*}

Each relation $\sigma _{i}^{m_{i}}=e$ in $P$ produces in $A^{\ast }$ a
relation of the form 
\begin{equation*}
\beta _{i}^{m_{i}}=\dprod\limits_{1\leq j\leq k}\beta _{j}^{xp_{ij}}
\end{equation*}%
where $p_{ij}$ are elements in the power series, as above; when written
additively $\beta _{i}^{m_{i}}$ has the form%
\begin{equation*}
m_{i}\beta _{i}=x\left( \sum_{1\leq j\leq k}p_{ij}\beta _{j}\right) \text{.}
\end{equation*}%
Let $F=\oplus _{1\leq i\leq k}\mathbb{Z}_{m}\left[ \left[ x\right] \right] 
\dot{\beta}_{i}$ be a free $\mathbb{Z}_{m}\left[ \left[ x\right] \right] $
module of rank $k$. Define the $\mathbb{Z}_{m}\left[ \left[ x\right] \right] 
$-homomorphism%
\begin{equation*}
\phi :\sum_{1\leq i\leq k}\mathbb{Z}_{m}\left[ \left[ x\right] \right] \dot{%
\beta}_{i}\rightarrow A^{\ast },\sum_{1\leq i\leq k}p_{i}\dot{\beta}%
_{i}\rightarrow \dprod\limits_{1\leq i\leq k}\beta _{i}^{p_{i}}
\end{equation*}%
and let $R$ be the kernel $\phi $. Define $J$ to be the $\mathbb{Z}_{m}\left[
\left[ x\right] \right] $-submodule of $R$ generated by 
\begin{equation*}
\dot{r}_{i}=m_{i}\dot{\beta}_{i}-x\left( \sum_{1\leq j\leq k}p_{ij}\dot{\beta%
}_{j}\right) \text{ }\left( 1\leq i\leq k\right) \text{.}
\end{equation*}%
We will show that $J=R$. So let $\nu \in R$ and write $\nu =\sum_{1\leq
i\leq k}\nu _{i}\dot{\beta}_{i}$ where%
\begin{eqnarray*}
\nu _{i} &=&\sum_{j\geq 0}\nu _{ij}x^{j}, \\
\nu _{ij} &=&\nu _{ij,0}+mw_{ij}\in \mathbb{Z}_{m}\text{.}
\end{eqnarray*}%
Then, $m_{i}|\nu _{i0,0}$, $\nu _{i0,0}=m_{i}\nu _{i0,0}^{\prime }$; factor $%
m=m_{i}m_{i}^{\prime }$. Therefore, 
\begin{eqnarray*}
\nu _{i} &=&\nu _{i0}+\left( \sum_{j\geq 1}\nu _{ij}x^{j-1}\right) x, \\
\nu _{i0} &=&m_{i}\nu _{i0,0}^{\prime }+mw_{i0}=\left( \nu _{i0,0}^{\prime
}+m_{i}^{\prime }w_{i0}\right) m_{i}, \\
\nu _{i}\dot{\beta}_{i} &=&\left( \nu _{i0,0}^{\prime }+m_{i}^{\prime
}w_{i0}\right) \left( m_{i}\dot{\beta}_{i}\right) +\left( \sum_{j\geq 1}\nu
_{ij}x^{j-1}\right) x\dot{\beta}_{i}, \\
&\equiv &\left( \nu _{i0,0}^{\prime }+m_{i}^{\prime }w_{i0}\right) \left(
x\sum_{1\leq j\leq k}p_{ij}\dot{\beta}_{j}\right) +\left( \sum_{j\geq 1}\nu
_{ij}x^{j-1}\right) x\dot{\beta}_{i}\text{ modulo }J\text{.}
\end{eqnarray*}%
Hence,

\begin{eqnarray*}
\nu &=&\sum_{1\leq i\leq k}\nu _{i}\dot{\beta}_{i}\in x\mu +J,\text{ } \\
\mu &=&\sum_{1\leq i\leq k}\mu _{i}\dot{\beta}_{i}\in R\text{.}
\end{eqnarray*}%
Thence, by repeating the argument, we obtain%
\begin{eqnarray*}
\nu &\in &\left( \dbigcap\limits_{i\geq 1}x^{i}R\right) +J=J\text{,} \\
J &=&R\text{. }
\end{eqnarray*}

On re-writing the relations $m_{i}\beta _{i}=\sum_{1\leq j\leq
k}p_{ij}x\beta _{j}$ in the form 
\begin{equation*}
p_{i1}x\beta _{1}+...+\left( p_{ii}x-m_{i}\right) \beta
_{i}+...+p_{kk}x\beta _{k}=0
\end{equation*}%
we see that the $k\times k$ matrix of coefficients of these equations has
determinant $r=m-xq$ for some $q\in \mathbb{Z}_{m}\left[ \left[ x\right] %
\right] $ and thus $r$ annuls $A^{\ast }$.

The last assertion of the theorem follows by using $r=m-qx\in R$ to reduce
the coefficients modulo $m$.
\end{proof}

\section{The $m$-congruence property}

A group $G$ of automorphisms of the $m$-ary tree is said to satisfy the $m$-%
\textit{congruence property} provided given $m^{i}$ there exists $l\left(
i\right) \geq 1$ such that $Stab_{G}(l\left( i\right) )\leq G^{m^{i}}$ for
all $i$; in which case the topology on $G$ inherited from $\mathcal{A}\left(
Y\right) $ is equal to the pro-$m$ topology. Since when $A^{\ast }$ is
written additively, we have $Stab_{G}(l\left( i\right) )=x^{l\left( i\right)
}A^{\ast }$, the $m$-congruence property reads $x^{l\left( i\right) }A^{\ast
}\leq m^{i}A^{\ast }$.

\begin{theorem}
Let $r=m-qx^{j}$ $\in \mathbb{Z}_{m}\left[ \left[ x\right] \right] $ with $%
q\in \mathbb{Z}_{m}\left[ \left[ x\right] \right] $ and $j\geq 1$. Let $S$
be quotient ring $\frac{\mathbb{Z}_{m}\left[ \left[ x\right] \right] }{%
\left( r\right) }$. Suppose $S$ is torsion-free. Then, $S$ is a finitely
generated pro-$m$ group.
\end{theorem}

\begin{proof}
From the decomposition $\mathbb{Z}_{m}\left[ \left[ x\right] \right] =\oplus
_{1\leq i\leq s}\varepsilon _{i}\mathbb{Z}_{p_{i}^{k_{i}}}\left[ \left[ x%
\right] \right] $ corresponding to the prime decomposition $m$ $%
=\dprod\limits_{1\leq i\leq s}p_{i}^{k_{i}}$, we obtain 
\begin{eqnarray*}
r &=&\sum_{1\leq i\leq s}r_{i}, \\
r_{i} &=&\varepsilon _{i}r=p_{i}^{k_{i}}-q_{i}\left( x\right) x^{j}, \\
S &=&\sum_{1\leq i\leq s}S_{i},\text{ \ }S_{i}=\frac{\mathbb{Z}%
_{p_{i}^{k_{i}}}\left[ \left[ x\right] \right] }{\left( r_{i}\right) }
\end{eqnarray*}%
where each $S_{i}$ is torsion-free. Thus, it is sufficient to address the
case where $m$ is a prime power $p^{k}$.

(1) First, we show that $S$ is a pro-$m$ group.

So, let $r=p^{k}-qx^{j}$ and decompose $q=q\left( x\right) =s\left( x\right)
+p.t\left( x\right) $ where each non-zero coefficient of $s\left( x\right) $
is an integer relatively prime to $p$. If $s\left( x\right) =0$ then $%
q\left( x\right) =p.t\left( x\right) $ and 
\begin{eqnarray*}
r &=&p^{k}-q\left( x\right) x^{j}=p^{k}-p.t\left( x\right) x^{j} \\
&=&p\left( p^{k-1}-t\left( x\right) x^{j}\right) \in \left( r\right) \text{;}
\end{eqnarray*}%
but as by hypothesis $S$ is torsion free, we have $p^{k-1}-t\left( x\right)
x^{j}\in \left( r\right) $ which is not possible.

Write $s\left( x\right) =x^{l}u\left( x\right) $ where $l\geq 0$ and where $%
u\left( x\right) $ is invertible in $\mathbb{Z}_{m}\left[ \left[ x\right] %
\right] $ with inverse $u^{\prime }\left( x\right) $. Then, $q\left(
x\right) =x^{l}u\left( x\right) +p.t\left( x\right) $ and 
\begin{eqnarray*}
r &=&p^{k}-\left( x^{l}u\left( x\right) x^{j}+p.t\left( x\right) x^{j}\right)
\\
&=&p\left( p^{k-1}-t\left( x\right) x^{j}\right) -x^{j+l}u\left( x\right) 
\text{.}
\end{eqnarray*}%
Therefore, on multiplying by $u^{\prime }\left( x\right) $, the inverse of $%
u\left( x\right) $, we obtain 
\begin{equation*}
p\left( p^{k-1}-t\left( x\right) x^{j}\right) u^{\prime }\left( x\right)
\equiv x^{j+l}\text{ modulo }r\text{.}
\end{equation*}%
It follows that 
\begin{equation*}
x^{j+l}S\leq pS,\text{ \ }x^{n\left( j+l\right) }S\leq p^{n}S\text{.}
\end{equation*}

(2) Now we show that $S$ is finitely generated as a $\mathbb{Z}_{m}$-module.

By the previous step there exist $l\geq 1$ and $v\left( x\right) \in \mathbb{%
Z}\left[ \left[ x\right] \right] $ such that 
\begin{equation*}
x^{l}\equiv mv\left( x\right) \text{ mod }r\text{.}
\end{equation*}%
Decompose $v\left( x\right) =v_{1}\left( x\right) +v_{2}\left( x\right)
x^{l} $ where the degree of $v_{1}\left( x\right) $ is less that $l$. Then,
we deduce modulo $r$: 
\begin{eqnarray*}
v\left( x\right) &\equiv &v_{1}\left( x\right) +v_{2}\left( x\right)
mv\left( x\right) , \\
v_{2}\left( x\right) v\left( x\right) &\equiv &w\left( x\right) \in \mathbb{Z%
}\left[ \left[ x\right] \right] , \\
w\left( x\right) &=&w_{1}\left( x\right) +w_{2}\left( x\right) x^{l}, \\
v\left( x\right) &\equiv &v_{1}\left( x\right) +mw\left( x\right) \\
&\equiv &v_{1}\left( x\right) +mw_{1}\left( x\right) +mw_{2}\left( x\right)
x^{l} \\
&&... \\
v\left( x\right) &\equiv &a_{0}+a_{1}x+...+a_{1}x^{l-1}\text{, }a_{i}\in 
\mathbb{Z}_{m}\text{.}
\end{eqnarray*}

We have shown that $S$ is generated by $1,x,...,x^{l-1}$ as a pro-$m$ group.
\end{proof}

\begin{corollary}
Let $A$ be an abelian transitive state-closed group of degree $m$. Suppose
the group $A^{\ast }$ is torsion-free. Then $A^{\ast }$ is a finitely
generated pro-$m$ group.
\end{corollary}

\begin{proof}
With previous notation, the group $A^{\ast }$ is a $\mathbb{Z}_{m}\left[ %
\left[ x\right] \right] $-module generated by 
\begin{equation*}
\left\{ \beta _{i}=\beta \left( \sigma _{i}\right) \mid 1\leq i\leq k\right\}
\end{equation*}%
and is annihilated by $r=m-qx^{j}$ $\in \mathbb{Z}_{m}\left[ \left[ x\right] %
\right] $ for some $q\in \mathbb{Z}_{m}\left[ \left[ x\right] \right] $ and $%
j\geq 1$.

It follows that $A^{\ast }$ is an $S$-module where $S=\frac{\mathbb{Z}_{m}%
\left[ \left[ x\right] \right] }{\left( r\right) }$. Since $S$ satisfies the 
$m$-congruence property, it follows that $A^{\ast }$ is a pro-$m$ group.

That $A^{\ast }$ is a finitely generated $\mathbb{Z}_{m}$-module, is a
consequence of $S$ being a finitely generated $\mathbb{Z}_{m}$-module.
\end{proof}

\section{Torsion in state-closed abelian groups}

\subsection{Preliminaries on virtual endomorphisms of groups}

Let $G$ be a transitive state-closed subgroup of $\mathcal{A}\left( Y\right) 
$ where $Y=\left\{ 1,2,...,m\right\} $. Then $\left[ G:Fix_{G}\left(
1\right) \right] =m$ and the projection on the $1$st coordinate of $%
Fix_{G}\left( 1\right) $ produces a subgroup of $G$; that is, $\pi
_{1}:Fix_{G}\left( 1\right) \rightarrow G$ is a virtual endomorphism of $G$.
This notion has proven to be effective in studying state-closed groups. We
give a quick review below.

Let $G$ be a group with a subgroup $H$ of finite index $m$ and a
homomorphism $f:H\rightarrow G$. A subgroup $U$ of $G$ is \textit{%
semi-invariant} under the action of $f$ provided $\left( U\cap H\right)
^{f}\leq U$. If $U\leq H$ and $U^{f}\leq U$ then $U$ is $f$-\textit{invariant%
}.

The largest subgroup $K$ of $H$, which is normal in $G$ and is $f$-invariant
is called the $f$-core$\left( H\right) $. If the $f$-core$\left( H\right) $
is trivial then $f$ and the triple $\left( G,H,f\right) $ are said to be a 
\textit{simple.}

Given a triple $\left( G,H,f\right) $ and a right transversal $L=\left\{
x_{1},x_{2},...,x_{m}\right\} $ of $H$ in $G$, the permutational
representation $\pi :G\rightarrow Perm\left( 1,2,...,m\right) $ is $g^{\pi
}:i\rightarrow j$ which is induced from the right multiplication $%
Hx_{i}g=Hx_{j}$. We produce recursively a representation $\varphi
:G\rightarrow $ $\mathcal{A}\left( m\right) $ as follows: 
\begin{equation*}
g^{\varphi }=\left( \left( x_{i}g.\left( x_{\left( i\right) g^{\pi }}\right)
^{-1}\right) ^{f\varphi }\right) _{1\leq i\leq m}g^{\pi }\text{.}
\end{equation*}%
One further expansion of $g^{\varphi }$ is 
\begin{eqnarray*}
g^{\varphi } &=&\left( \left( \left( x_{j}g_{i}.x_{\left( j\right)
g_{i}^{\pi }}^{-1}\right) ^{f\varphi }\right) _{1\leq j\leq m}g_{i}^{\pi
}\right) _{1\leq i\leq m}g^{\pi }, \\
&=&\left( \left( \left( x_{j}g_{i}.x_{\left( j\right) g_{i}^{\pi
}}^{-1}\right) ^{f\varphi }\right) _{1\leq j\leq m}\right) _{1\leq i\leq
m}\left( g_{i}^{\pi }\right) _{1\leq i\leq m}\text{ }g^{\pi }
\end{eqnarray*}%
where $g_{i}=\left( x_{i}g.x_{\left( i\right) g^{\pi }}^{-1}\right) ^{f}$

The kernel of $\varphi $ is precisely the $f$-core$\left( H\right) $, $%
G^{\varphi }$ is state-closed and $H^{\varphi }=Fix_{G^{\varphi }}(1)$.

\subsubsection{Changing transversals}

We will show below that changing the transversal of $H$ in $G$ produces
another representation of $G$, conjugate to the original one by an explicit
automorphism of the $m$-ary tree.

\begin{proposition}
Let $\left( G,H,f\right) $ be a triple and 
\begin{equation*}
L=\left\{ x_{1},x_{2},...,x_{m}\right\} ,L^{\prime }=\left\{ x_{1}^{\prime
}=h_{1}x_{1},x_{2}^{\prime }=h_{2}x_{2},...,x_{m}^{\prime
}=h_{m}x_{m}\right\}
\end{equation*}%
right transversals of $H$ in $G$ where $h_{i}\in H$. Let $\varphi =\varphi
_{x_{i}},\varphi ^{\prime }=\varphi _{h_{i}x_{i}}:G\rightarrow $ $\mathcal{A}%
\left( m\right) $ be the corresponding tree representations and define the
following elements of $\mathcal{A}\left( m\right) $, 
\begin{eqnarray*}
\gamma &=&\gamma _{h_{i},\varphi ^{\prime }}=\left( \left( h_{i}\right)
^{f\varphi ^{\prime }}\right) _{1\leq i\leq m}, \\
\lambda &=&\lambda _{h_{i},\varphi ^{\prime }}=\gamma \gamma ^{\left(
1\right) }...\gamma ^{\left( n\right) }...\text{.}
\end{eqnarray*}%
Then,%
\begin{equation*}
\varphi _{h_{i}x_{i}}=\varphi _{x_{i}}\left( \lambda _{h_{i}^{-1},\varphi
_{x_{i}}}\right) \text{.}
\end{equation*}
\end{proposition}

\begin{proof}
The representations $\varphi ,\varphi ^{\prime }:G\rightarrow $ $\mathcal{A}%
\left( m\right) $ are defined by 
\begin{eqnarray*}
g^{\varphi } &=&\left( \left( x_{i}g.\left( x_{\left( i\right) g^{\pi
}}\right) ^{-1}\right) ^{f\varphi }\right) _{1\leq i\leq m}g^{\pi }\text{,}
\\
g^{\varphi ^{\prime }} &=&\left( \left( x_{i}^{\prime }g.\left( x_{\left(
i\right) g^{\pi }}^{\prime }\right) ^{-1}\right) ^{f\varphi ^{\prime
}}\right) _{1\leq i\leq m}g^{\pi }\text{.}
\end{eqnarray*}%
The relationship between $\varphi ^{\prime }$ and $\varphi $ is established
as follows, 
\begin{eqnarray*}
g^{\varphi ^{\prime }} &=&\left( \left( h_{i}x_{i}g.\left( h_{\left(
i\right) g^{\pi }}x_{\left( i\right) g^{\pi }}\right) ^{-1}\right)
^{f\varphi ^{\prime }}\right) _{1\leq i\leq m}g^{\pi } \\
&=&\left( \left( h_{i}\left( x_{i}g.x_{\left( i\right) g^{\pi }}^{-1}\right)
h_{\left( i\right) g^{\pi }}^{-1}\right) ^{f\varphi ^{\prime }}\right)
_{1\leq i\leq m}g^{\pi } \\
&=&\left( \left( h_{i}\right) ^{f\varphi ^{\prime }}\right) _{1\leq i\leq
m}.\left( \left( x_{i}g.x_{\left( i\right) g^{\pi }}^{-1}\right) ^{f\varphi
^{\prime }}\right) _{1\leq i\leq m}.\left( \left( h_{\left( i\right) g^{\pi
}}^{-1}\right) ^{f\varphi ^{\prime }}\right) _{1\leq i\leq m}g^{\pi } \\
&=&\left( \left( h_{i}\right) ^{f\varphi ^{\prime }}\right) _{1\leq i\leq
m}.\left( \left( x_{i}g.x_{\left( i\right) g^{\pi }}^{-1}\right) ^{f\varphi
^{\prime }}\right) _{1\leq i\leq m}g^{\pi }.\left( \left( h_{i}\right)
^{f\varphi ^{\prime }}\right) _{1\leq i\leq m}^{-1}\text{.}
\end{eqnarray*}%
Therefore 
\begin{equation*}
g^{\varphi ^{\prime }}=\gamma .\left( \left( x_{i}g.x_{\left( i\right)
g^{\pi }}^{-1}\right) ^{f\varphi ^{\prime }}\right) _{1\leq i\leq m}g^{\pi
}.\gamma ^{-1}
\end{equation*}%
where $\gamma =\left( \left( h_{i}\right) ^{f\varphi ^{\prime }}\right)
_{1\leq i\leq m}$ is independent of $g$. Repeating this development for each 
$g_{i}=\left( x_{i}g.x_{\left( i\right) g^{\pi }}^{-1}\right) ^{f}$, we find
that 
\begin{equation*}
g^{\varphi ^{\prime }}=\gamma \gamma ^{\left( 1\right) }.\left( \left(
\left( x_{j}g_{i}.x_{\left( j\right) g_{i}^{\pi }}^{-1}\right) ^{f\varphi
^{\prime }}\right) _{1\leq j\leq m}g_{i}^{\pi }\right) _{1\leq i\leq
m}g^{\pi }.\gamma ^{-\left( 1\right) }\gamma ^{-1}\text{.}
\end{equation*}%
Thus in the limit, we obtain $\lambda =\gamma \gamma ^{\left( 1\right)
}...\gamma ^{\left( n\right) }...$ such that%
\begin{eqnarray*}
g^{\varphi ^{\prime }} &=&\lambda g^{\varphi }\lambda ^{-1}\text{ for all }%
g\in G, \\
\varphi &=&\varphi ^{\prime }\lambda \text{.}
\end{eqnarray*}

Introducing the explicit dependence of $\varphi ,\varphi ^{\prime },\lambda $
on the transversals, the previous equation becomes 
\begin{equation*}
\varphi _{x_{i}}=\left( \varphi _{h_{i}x_{i}}\right) \left( \lambda
_{h_{i},\varphi _{h_{i}x_{i}}}\right) \text{.}
\end{equation*}%
On replacing $h_{i}$ by $h_{i}^{-1}$ and on denoting $h_{i}^{-1}x_{i}$ by $%
x_{i}^{\prime }$, we obtain

\begin{equation*}
\varphi _{h_{i}x_{i}^{\prime }}=\left( \varphi _{x_{i}^{\prime }}\right)
\left( \lambda _{h_{i}^{-1},\varphi _{x_{i}^{\prime }}}\right) \text{.}
\end{equation*}
\end{proof}

\begin{example}
Let $G=C=\left\langle a\right\rangle $ be the infinite cyclic group, let $%
H=\left\langle a^{2}\right\rangle $ and let $f:H\rightarrow G$ be defined by 
$a^{2}\rightarrow a$. Given $l,k\geq 0$, then on choosing the transversal $%
L_{k,l}=\left\{ a^{2k},a^{2l+1}\right\} $ for $H$ in $G$, we obtain the
representation $\varphi _{k,l}:G\rightarrow \mathcal{A}\left( m\right) $
where $\varphi _{k,l}:a\rightarrow \alpha =\left( \alpha ^{k-l},\alpha
^{-k+l+1}\right) \sigma $.
\end{example}

\subsubsection{Subtriples, Quotient triples}

Given a triple $\left( G,H,f\right) $ and given subgroups $V\leq G,U\leq
H\cap V$ such that $\left( U\right) ^{f}\leq V$, we call $\left(
V,U,f|_{U}\right) $ a \textit{sub-triple }of $G$\textit{. }If $N$ is a
normal semi-invariant subgroup of $G$ then $\overline{f}:\frac{HN}{N}%
\rightarrow \frac{G}{N}$ given by $\overline{f}:Nh\rightarrow Nh^{f}$ is
well-defined and $\left( \frac{G}{N},\frac{HN}{N},\overline{f}\right) $ is a 
\textit{quotient triple}.

Let $\left( G,H,f\right) $ be a simple triple where $G$ is abelian and $%
\left[ G:H\right] =m$. Then, any sub-triple of $G$ is simple. Let $%
T=tor\left( G\right) $ denote the torsion subgroup of $G$ and for $l\geq 1$
define $G\left( l\right) =\left\{ g\in T\mid o\left( g\right) |l\right\}
,H\left( l\right) =G\left( l\right) \cap H$. Then, clearly, $f:tor\left(
H\right) \rightarrow tor\left( G\right) $ and $f:H\left( l\right)
\rightarrow G\left( l\right) $. Therefore, $tor\left( G\right) $ and $%
G\left( l\right) $ are semi-invariant and $\left( tor\left( G\right)
,tor\left( H\right) ,f|_{tor\left( H\right) }\right) $ and $\left( G\left(
l\right) ,H\left( l\right) ,f|_{H\left( l\right) }\right) $ are simple 
\textit{sub-triples}.

\begin{lemma}
Let $\left( G,H,f\right) $ be a simple triple. The triple $\left( \frac{G}{%
G\left( l\right) },\frac{HG\left( l\right) }{G\left( l\right) },\overline{f}%
\right) $ is also simple.
\end{lemma}

\begin{proof}
For suppose $K\leq H$ is such that $G\left( l\right) K^{f}\leq G\left(
l\right) K$. Then 
\begin{equation*}
\left( G\left( l\right) K^{f}\right) ^{l}=\left( K^{f}\right) ^{l}=\left(
K^{l}\right) ^{f}\leq \left( G\left( l\right) K\right) ^{l}=\left( K\right)
^{l}\text{;}
\end{equation*}%
that is, $K^{l}$ is $f$-invariant. Since $f$ is simple, $K^{l}=\left\{
e\right\} $ and so, $K\leq G\left( l\right) $.
\end{proof}

\subsection{The torsion subgroup}

\begin{proposition}
Let $A$ be transitive state-closed abelian group of degree $m$. Then $%
tor\left( A\right) $ has finite exponent and is therefore a direct summand
of $A$.
\end{proposition}

\begin{proof}
Let $T=tor\left( A\right) ,$ $A_{1}=Stab_{A}\left( 1\right) $, $T_{1}=T\cap
A_{1}$and $\left[ T:T_{1}\right] =m^{\prime }$. Then, the projection on the $%
1$st coordinate of $T_{1}$ is a subgroup of $T$ and the triple $\left(
T,T_{1},\pi _{1}|_{T_{1}}\right) $ is simple of degree $m^{\prime }|m$; let $%
m=m^{\prime }m^{\prime \prime }$. Hence, in this representation, $T$ is a
torsion transitive state-closed subgroup of $\mathcal{A}_{m^{\prime }}$, the
automorphism group of the tree $\mathcal{T}_{m^{\prime }}$.

Fixing this last representation of $T$, let $Q=P\left( T\right) $ and let $%
\sigma _{i}$ $\left( 1\leq i\leq k\right) $ be a minimal set of generators
of $Q$ and as before, let $\beta _{i}=\beta \left( \sigma _{i}\right) \in T$
be such that $\sigma \left( \beta _{i}\right) =\sigma _{i}$. Let $r$ be the
maximum order of the elements $\beta _{1},...,\beta _{k}$. As any $\alpha
\in T$ can be written in the form%
\begin{equation*}
\alpha =\dprod\limits_{1\leq i\leq k}\beta _{i}^{r_{i1}}\left( \beta
_{i}^{r_{i2}}\right) ^{\left( l_{2}\right) }...\left( \beta
_{i}^{r_{ij}}\right) ^{\left( l_{j}\right) }...
\end{equation*}%
it follows that $\alpha ^{r}=e$.

Since $T$ has finite exponent, it is a pure bounded subgroup of $A$ and
therefore it is a direct summand of $A$ (\cite{Rob}, Th. 4.3.8).
\end{proof}

We recall a classic example of an abelian group $G$ which does not split
over its torsion subgroup (see \cite{Rob}, page 108).

\begin{example}
Let $G$ be the direct product of groups $\dprod\limits_{i\geq 1}C_{i}$ where 
$C_{i}=\left\langle c_{i}\right\rangle $ is cyclic of order $p^{i}$ and let $%
H$ be the direct sum $\dsum\limits_{i\geq 1}C_{i}$. Then $H\leq tor\left(
G\right) =\dbigcup\limits_{l\geq 1}G\left( p^{l}\right) $. Moreover, $H$ is
a basic subgroup of $G$ and in particular, $\frac{G}{H}$ is $p$-divisible.
This observation leads directly to a proof that $G$ does not split over $%
tor\left( G\right) $.
\end{example}

The proof of the previous proposition did not establish the exponent of $%
tor\left( A\right) $. This we do in the next two lemmas.

\begin{lemma}
Let $m$ be a prime number and $A$ an abelian transitive state-closed \textit{%
torsion} group of degree $m$. Then, $A$ is conjugate by a tree automorphism
to a subgroup of the diagonal-topological closure of $\left\langle \sigma
\right\rangle $ and so has exponent $m$.
\end{lemma}

\begin{proof}
We observe that, $A\left( m\right) $ is not contained in $%
A_{1}=Stab_{A}\left( 1\right) $. For otherwise, $A\left( m\right) $ would be
invariant under the projection on the $1$st coordinate. Choose $a\in
A\backslash A_{1}$ of order $m$. Therefore, $A=A_{1}\left\langle
a\right\rangle $. On choosing $\left\{ a^{i}\mid 0\leq i\leq m-1\right\} $
as a transversal of $A_{1}$ in $A$, the image of $a$ acquires the form $%
\sigma =(1,...,m)$ in this tree representation of $A$. Thus, we may suppose
by Proposition 3 that $\sigma \in A$. Therefore, $\widetilde{A}$ contains
the subgroup $\widetilde{\left\langle \sigma \right\rangle }=\left\langle
\sigma ^{\left( i\right) }\mid i\geq 0\right\rangle $. By Proposition 2, we
have $C_{\mathcal{A}}\widetilde{\left\langle \sigma \right\rangle }%
=\left\langle \sigma \right\rangle ^{\ast }$ and thus, $A\leq C_{\mathcal{A}%
}\left( A\right) \leq \left\langle \sigma \right\rangle ^{\ast }$.
\end{proof}

\begin{lemma}
Suppose $A$ is an abelian transitive state-closed torsion group of degree $m$%
. Then the exponent of $A$ is equal to the exponent of $P\left( A\right) $.
\end{lemma}

\begin{proof}
By induction on $|P\left( A\right) |=m$. The exponent of $A$ is a multiple
of the exponent of $P\left( A\right) $. By the previous lemma, we may assume 
$m$ to be composite. Let $p$ be a prime divisor of $m$ and $A\left( p\right)
=\left\{ a\in A\mid a^{p}=e\right\} $. Then, $A\left( p\right) $ is a
nontrivial subgroup and $P\left( A\left( p\right) \right) \leq \left\{
\sigma \in P\mid \sigma ^{p}=e\right\} $. By Lemma 2, $\left( \frac{A}{%
A\left( p\right) },\frac{A_{1}A\left( p\right) }{A\left( p\right) },%
\overline{\pi _{0}}\right) $ is simple; also, $P\left( \frac{A}{A\left(
p\right) }\right) =\frac{P\left( A\right) }{P\left( A\left( p\right) \right) 
}$. The proof follows by induction.
\end{proof}

\begin{theorem}
Suppose $A$ is an abelian transitive state-closed torsion group of degree $m$%
. Then, $A$ is conjugate to a subgroup of the topological closure of%
\begin{equation*}
\widetilde{P\left( A\right) }=\left\langle \sigma ^{\left( i\right) }\mid
\sigma \in P\left( A\right) \text{, }i\geq 0\right\rangle \text{.}
\end{equation*}
\end{theorem}

\begin{proof}
Let $P=P\left( A\right) $ have exponent $r$ and let $B$ be a maximal
homogeneous subgroup of $P$ of exponent $r$ (that is, $B$ is a direct sum of
cyclic groups of order $r)$, minimally generated by $\left\{ \sigma _{i}\mid
\left( 1\leq i\leq s\right) \right\} $. Choose for each $\sigma _{i}$ an
element $\beta _{i}=\beta \left( \sigma _{i}\right) \in A$ and let $\dot{B}%
=\left\langle \beta _{i}\mid \left( 1\leq i\leq s\right) \right\rangle $.
Then, as the order of each $\beta _{i}$ is a multiple of $r$, while the
exponent of $A$ is $r$, we conclude from the previous lemma that $o\left(
\beta _{i}\right) =o\left( \sigma _{i}\right) =r$ for $\left( 1\leq i\leq
s\right) $. Since $\beta _{i}\rightarrow \sigma _{i}$ defines a projection
of $\dot{B}$ onto $B$ we conclude that $\dot{B}\cong B$ and $\dot{B}\cap
A_{1}=\left\{ e\right\} $, where $A_{1}=Stab_{A}\left( 1\right) $.

Clearly $\dot{B}$ is a pure bounded subgroup and so it has a complement $L$
in $A$, which may be chosen to contain $A_{1}$. Choose a right transversal $%
W $ of $A_{1}$ in $L$. Then the set $W\dot{B}$ is a right transversal of $%
A_{1} $ in $A$. With respect to this transversal, the triple $\left(
A,A_{1},\pi _{1}\right) $ produces a transitive state-closed representation $%
\varphi $ where $\dot{B}^{\varphi }=B$. By Proposition 3, we may rewrite $%
A^{\varphi }$ as $A$. Then, the diagonal-topological closure $A^{\ast }$
contains $B^{\ast }$. Let $V$ be a complement of $B$ in $P$. Each $\alpha
\in A^{\ast }$ can be factored as $\alpha =\beta \gamma $ where $\beta \in
B^{\ast }$ and $\gamma $ is such that each of its states $\gamma _{u}$ have
activity $\sigma \left( \gamma _{u}\right) \in V$. Therefore, the set of
these $\gamma $'s is a group $\Gamma $ such that $\Gamma =\Gamma ^{\ast }$
and $A^{\ast }=\Gamma \oplus B^{\ast }$. Then, $\left( \Gamma ,\Gamma \cap
A_{1},\pi _{1}\right) $ is a simple triple with $P\left( \Gamma \right) $
having exponent smaller than $r$. The proof is finished by induction on the
exponent.
\end{proof}

The example below illustrates some of the ideas developed so far.

\begin{example}
Let $m=4,Y=\left\{ 1,2,3,4\right\} $and let $\sigma $ be the cycle $\left(
1,2,3,4\right) $. Furthermore, let $\alpha =\left( e,e,e,\alpha ^{2}\right)
\sigma \in \mathcal{A}\left( 4\right) $ and let $A=\left\langle \alpha
\right\rangle $. Then 
\begin{eqnarray*}
\alpha ^{2} &=&\left( \alpha ^{2},e,e,\alpha ^{2}\right) \left( 1,3\right)
\left( 2,4\right) , \\
\alpha ^{4} &=&\left( \alpha ^{2}\right) ^{\left( 1\right) }=\alpha
^{2x},\left( \alpha ^{2-x}\right) ^{2}=e\text{.}
\end{eqnarray*}%
Then $A$ is cyclic, torsion-free, transitive and state-closed, yet it is not
diagonally closed, as $\alpha ^{x}\not\in A$. Even though $A$ is
torsion-free, its diagonal closure $\widetilde{A}=\left\langle \alpha
^{x^{i}}\mid i\geq 0\right\rangle $ is not; for $\kappa =\alpha ^{2-x}$ has
order $2$. Let $K=\left\langle \kappa ^{x^{i}}\mid i\geq 0\right\rangle $.
Then, $K\leq tor\left( \widetilde{A}\right) $ and it is direct to check that 
$\widetilde{A}=\left\langle \alpha ,K\right\rangle $. Therefore, $%
K=tor\left( \widetilde{A}\right) $ and 
\begin{equation*}
\widetilde{A}=tor\left( \widetilde{A}\right) \oplus A\text{.}
\end{equation*}%
Let $Y_{1}=\left\{ 1,3\right\} ,Y_{2}=\left\{ 2,4\right\} $. Then $\left\{
Y_{1},Y_{2}\right\} $ is a complete block system for the action of $\alpha $
on $Y$. Also, $\alpha ^{2}$ induces the binary adding machine on both $%
\mathcal{T}\left( Y_{1}\right) $ and $\mathcal{T}\left( Y_{2}\right) $. The
topological closure $\overline{A}$ of $A$ is torsion-free and 
\begin{eqnarray*}
tor\left( A^{\ast }\right) &=&tor\left( \widetilde{A}\right) , \\
A^{\ast } &=&tor\left( A^{\ast }\right) \oplus \overline{A}\text{.}
\end{eqnarray*}%
Moreover, $tor\left( A^{\ast }\right) $ induces a faithful state-closed,
diagonally and topologically closed actions on the binary tree $\mathcal{T}%
\left( Y_{1}\right) $. Therefore, $tor\left( A^{\ast }\right) $ is
isomorphic to $\frac{\mathbb{Z}}{2\mathbb{Z}}\left[ \left[ x\right] \right] $%
. Furthermore, $\alpha $ is represented as the binary adding machine on $%
\mathcal{T}\left( \left\{ Y_{1},Y_{2}\right\} \right) $ and $\overline{A}$
is represented on this tree as the topological closure of the image of $A$.
\end{example}

\section{Cyclic $\mathbb{Z}_{m}\left[ \left[ x\right] \right] $-modules}

Cyclic automorphism groups $\left\langle \alpha \right\rangle $ of the tree,
for which their state-diagonal-topological closure is isomorphic to a cyclic 
$\mathbb{Z}_{m}$-module have the form 
\begin{equation*}
\alpha =\left( \alpha ^{q_{1}},...,\alpha ^{q_{m}}\right) \sigma
\end{equation*}%
where $q_{i}\in \mathbb{Z}_{m}\left[ \left[ x\right] \right] $ for $1\leq
i\leq m$; here%
\begin{eqnarray*}
q_{i} &=&\sum_{j\geq 0}q_{ij}x^{j}, \\
q_{ij} &=&\sum_{u\geq 0}q_{ij,u}m^{u}\in \mathbb{Z}_{m}\text{.}
\end{eqnarray*}

We prove

\begin{theorem}
(i) The expression 
\begin{equation*}
\alpha =\left( \alpha ^{q_{\mathbf{1}}},...,\alpha ^{q_{\mathbf{m}}}\right)
\sigma
\end{equation*}%
is a well-defined automorphism of the $m$-ary tree. (ii) Let $A$ be the
state closure of $\left\langle \alpha \right\rangle $. Then $A^{\ast }$is
abelian, isomorphic to the quotient ring $\frac{\mathbb{Z}_{m}\left[ \left[ x%
\right] \right] }{\left( r\right) }\ $where%
\begin{equation*}
r=m-qx\text{ and }q=q_{\mathbf{1}}+...+q_{\mathbf{m}}\text{.}
\end{equation*}
\end{theorem}

\begin{proof}
(1) Let $\sigma \left( l\right) $ denote the permutation induced by $\alpha $
on the $l$-th level. Then, the expression $\alpha =\left( \alpha ^{q_{%
\mathbf{1}}},...,\alpha ^{q_{\mathbf{m}}}\right) \sigma $ represents 
\begin{eqnarray*}
\sigma \left( 1\right) &=&\sigma , \\
\sigma \left( l\right) &=&\left( \sigma (l-1)^{\overline{q_{1}}},...,\sigma
(l-1)^{\overline{q_{m}}}\right) \sigma
\end{eqnarray*}%
where $\overline{q_{i}}=\overline{q_{i0}}+\overline{q_{i1}}x+...+\overline{%
q_{i(l-1)}}x^{l-1}$ and $\overline{q_{ij}}%
=q_{ij,0}+q_{ij,1}m+...+q_{ij,l-1}m^{l-1}$.

(2.1) The states of $\alpha $ are words in $\alpha ^{p}$ for $p\in \mathbb{Z}%
_{m}\left[ \left[ x\right] \right] $. Let $v=\alpha ^{l_{1}}...\alpha
^{l_{a}},w=\alpha ^{n_{1}}...\alpha ^{n_{b}}\in A^{\ast }$. Then clearly $%
\left[ v,w\right] \in Stab_{A}(1)$. We will prove that the entries of $\left[
v,w\right] $ are products of conjugates of words in elements of the form $%
\left[ \alpha ^{s},\alpha ^{t}\right] $ where $s,t\in \mathbb{Z}_{m}\left[ %
\left[ x\right] \right] $.

Clearly $\left[ v,w\right] $ can be developed into a word in conjugates of $%
\left[ \alpha ^{l_{i}},\alpha ^{n_{j}}\right] $.

Write $p=p_{0}+p^{\prime }x,n=n_{0}+n^{\prime }x$. We compute 
\begin{eqnarray*}
\left[ \alpha ^{p},\alpha ^{n}\right] &=&\left( \left[ \alpha
^{p_{0}},\alpha ^{n^{\prime }x}\right] \left[ \alpha ^{p_{0}},\alpha ^{n_{0}}%
\right] ^{\alpha ^{n^{\prime }x}}\right) ^{\alpha ^{p^{\prime }x}} \\
&&\left[ \alpha ^{p^{\prime }},\alpha ^{n^{\prime }}\right] ^{x}\left[
\alpha ^{p^{\prime }x},\alpha ^{n_{0}}\right] ^{\alpha ^{n^{\prime }x}} \\
&=&\left[ \alpha ^{p_{0}},\alpha ^{n^{\prime }x}\right] ^{\alpha ^{p^{\prime
}x}}\left[ \alpha ^{p^{\prime }},\alpha ^{n^{\prime }}\right] ^{x}\left[
\alpha ^{p^{\prime }x},\alpha ^{n_{0}}\right] ^{\alpha ^{n^{\prime }x}}\text{%
.}
\end{eqnarray*}%
Therefore, we have to check $\left[ \alpha ^{\xi },\alpha ^{nx}\right] $
where $\xi \in \mathbb{Z}_{m},n\in \mathbb{Z}_{m}\left[ \left[ x\right] %
\right] $. Write $\xi =\xi _{0}+m\xi ^{\prime }$. Then, 
\begin{eqnarray*}
\left[ \alpha ^{\xi },\alpha ^{nx}\right] &=&\left[ \alpha ^{\xi _{0}+m\xi
^{\prime }},\alpha ^{nx}\right] \\
&=&\left[ \alpha ^{\xi _{0}},\alpha ^{nx}\right] ^{\alpha ^{m\xi ^{\prime }}}%
\left[ \alpha ^{m\xi ^{\prime }},\alpha ^{nx}\right] \text{.}
\end{eqnarray*}%
Now, 
\begin{equation*}
\alpha ^{\xi _{0}}=\left( v_{\mathbf{1}},v_{\mathbf{2}},...,v_{\mathbf{m}%
}\right) \sigma ^{\xi _{0}},
\end{equation*}%
where $v_{\mathbf{i}}$ are words in $\alpha ^{q_{\mathbf{1}}},...,\alpha
^{q_{\mathbf{m}}}$ and%
\begin{equation*}
\alpha ^{m}=\left( \alpha ^{q_{\mathbf{1}}}...\alpha ^{q_{\mathbf{m}%
}},\alpha ^{q_{\mathbf{2}}}...\alpha ^{q_{\mathbf{m}}}\alpha ^{q_{\mathbf{1}%
}},...,\alpha ^{q_{\mathbf{m}}}\alpha ^{q_{\mathbf{1}}}...\alpha ^{q_{%
\mathbf{m-1}}}\right) \text{.}
\end{equation*}%
Therefore, 
\begin{equation*}
\left[ \alpha ^{\xi _{0}},\alpha ^{nx}\right] =\left( \left[ v_{\mathbf{1}%
},\alpha ^{n}\right] ,...,\left[ v_{\mathbf{m}},\alpha ^{n}\right] \right)
\end{equation*}%
and similarly,%
\begin{equation*}
\left[ \alpha ^{m\xi ^{\prime }},\alpha ^{nx}\right] =\left( \left[ \left(
\alpha ^{q_{\mathbf{1}}}...\alpha ^{q_{\mathbf{m}}}\right) ^{\xi ^{\prime
}},\alpha ^{n}\right] ,...,\left[ \left( \alpha ^{q_{\mathbf{m}}}\alpha ^{q_{%
\mathbf{1}}}...\alpha ^{q_{\mathbf{m-1}}}\right) ^{\xi ^{\prime }},\alpha
^{n}\right] \right) \text{.}
\end{equation*}%
Now we write $\beta =\alpha ^{q_{\mathbf{1}}}...\alpha ^{q_{\mathbf{m}}}$.
Then $\left[ \beta ^{\xi ^{\prime }},\alpha ^{n}\right] $ can be developed
further as asserted. The same applies to the other entries.

(2.2) First, clearly $r\alpha =0$. Now let $u=u\left( x\right) $ annul $%
\alpha $; write $u=u_{0}+u^{\prime }x$ where $u_{0}=u\left( 0\right) $. Then 
$m|u_{0}$ and so,

\begin{eqnarray*}
u &=&m\frac{u_{0}}{m}+u^{\prime }x=\left( xq\right) \frac{u_{0}}{m}%
+u^{\prime }x+vr \\
&=&xw_{1}+vr
\end{eqnarray*}%
for some $v=v\left( x\right) $ and $w_{1}=q\frac{u_{0}}{m}+u^{\prime }$.
Then, $xw_{1}$ annuls $\alpha $ and so does $w_{1}$. On repeating, we find $%
w_{i}$ such that $u\equiv x^{i}w_{i}$ modulo $r$ and $w_{i}$ annuls $\alpha $
for all $i\geq 1$.

In other words, $u\in \cap _{n\geq 1}\left( x\mathbb{Z}\right) ^{n}+\left(
r\right) =\left( r\right) $.
\end{proof}

\subsubsection{The group $\mathit{D}_{m}\left( j\right) $}

Recall $\alpha =\left( e,..,e,\alpha ^{x^{j-1}}\right) \sigma \in \mathcal{A}%
_{m}$. Then $\alpha ^{m}=\alpha ^{x^{j}}$; that is, $\alpha ^{r}=e$ where $%
r=m-x^{j}$. The states of $\alpha $ are $\alpha ,\alpha ^{x},...,\alpha
^{x^{j-1}}$ and 
\begin{equation*}
\mathit{D}_{m}\left( j\right) =\left\langle \alpha ,\alpha ^{x},...,\alpha
^{x^{j-1}}\right\rangle \text{;}
\end{equation*}%
therefore $\mathit{D}_{m}\left( j\right) $ is diagonally closed. The
topological closure $\overline{\mathit{D}_{m}\left( j\right) }$ is
isomorphic to the quotient ring $S=\frac{\mathbb{Z}_{m}\left[ \left[ x\right]
\right] }{\left( r\right) }$ which is clearly a free $\mathbb{Z}_{m}$-module
of rank $j$.

\subsection{The case $P\left( A\right) $ cyclic of prime order}

\begin{theorem}
Let $m$ be a prime number. Let $A$ be a torsion-free abelian transitive
state-closed subgroup of $\mathcal{A}_{m}$. Let $\beta \in A\backslash
Stab_{A}\left( j\right) $. Then $A^{\ast }=$ $\left\langle \beta
\right\rangle ^{\ast }$ and is topologically finitely generated.
Furthermore, $A^{\ast }$ is conjugate to $\overline{\mathit{D}_{m}\left(
j\right) }$ for some $j\geq 1$.
\end{theorem}

The proof is developed in four steps.

\textbf{Step 1}. For $z\in A$, define $\zeta (z)=j$ such that $z^{m}\in
Stab\left( j\right) \backslash Stab\left( j+1\right) $. As $A$ is
torsion-free, $\zeta (z)$ is finite for all nontrivial $z$ and $z^{m}=\left(
v\right) ^{\left( j\right) },v\in A\backslash Stab_{A}\left( 1\right) $.

Choose $\beta =(\beta _{1},\beta _{2},...,\beta _{m})\sigma \in A\backslash
Stab_{A}\left( 1\right) $ having minimum $\zeta (\beta )=j$. If $z\in
Stab_{A}\left( 1\right) ,$ $z\not=e$, then there exists $l>0$ such that $%
z^{m}=(c)^{\left( l\right) }$ and $c\in A\backslash Stab_{A}\left( 1\right) $%
. Therefore, by minimality of $\beta $ we have $\zeta (c)\geq \zeta (\beta )$
and $\zeta (z)>\zeta (\beta )$.

\begin{lemma}
(Uniform gap) Let $z\in $ $Stab_{A}\left( 1\right) $. Then $\zeta (z\beta
)=\zeta (\beta )$.
\end{lemma}

\begin{proof}
First note that 
\begin{eqnarray*}
\beta ^{m} &=&(\beta _{1}\beta _{2}...\beta _{m})^{\left( 1\right) }, \\
\beta _{1}\beta _{2}...\beta _{m} &=&\left( \gamma \right) ^{\left(
j-1\right) }\text{, }\gamma \in A\backslash Stab_{A}\left( 1\right) \text{.}
\end{eqnarray*}%
We have $z=c^{\left( 1\right) }$ and $z\beta =(c\beta _{1},c\beta
_{2},...,c\beta _{m})\sigma ,(z\beta )^{m}=(u)^{\left( 1\right) }$ where $%
u=c^{m}\beta _{1}...\beta _{m}=c^{m}\left( \gamma \right) ^{\left(
j-1\right) }$. If $c\in A\backslash Stab_{A}\left( 1\right) $ then $\zeta
(c)=n\geq j,c^{m}\in Stab\left( n\right) \backslash Stab\left( n+1\right) $
and so, $u\in Stab_{A}\left( j-1\right) \backslash Stab_{A}\left( j\right) $%
. If $c\in Stab_{A}\left( 1\right) $ then $\zeta (c)>j$ and so, $c^{m}\in
Stab\left( k\right) $ where $k>j$ and again $u\in Stab\left( j-1\right)
\backslash Stab\left( j\right) $.
\end{proof}

\textbf{Step 2.} Note that 
\begin{eqnarray*}
\beta ^{m} &=&\left( \gamma \right) ^{\left( j\right) },\text{ \ }\gamma
^{m}=\left( \lambda \right) ^{\left( j\right) } \\
\beta ^{m^{2}} &=&\left( \lambda \right) ^{\left( 2j\right) }
\end{eqnarray*}%
where by the uniform gap lemma above, $\gamma ,\lambda \in A\backslash
Stab_{A}\left( 1\right) $. Therefore,\ on repeating this process, we find
that $\beta ^{m^{s}}$ induces $\sigma ^{\left( sj\right) }$ on the $\left(
sj\right) $th level of the tree for all $s\geq 0$. Now given a level $t\geq
0 $, on dividing $t$ by $j$, we get $t=sj+i$ with $0\leq i\leq j-1$ and then 
$\left( \beta ^{\left( i\right) }\right) ^{m^{s}}=\left( \beta
^{m^{s}}\right) ^{\left( i\right) }$ induces $\left( \sigma ^{\left(
sj\right) }\right) ^{\left( i\right) }=\sigma ^{\left( sj+i\right) }=\sigma
^{\left( t\right) }$ on the $t$-th level of the tree. It follows that the
group $A$ is a subgroup of the topological closure of $\left\langle \beta
,\beta ^{\left( 1\right) },....,\beta ^{\left( j-1\right) }\right\rangle $.

\textbf{Step 3.} We have for $\beta =(\beta _{1},\beta _{2},...,\beta
_{m})\sigma $, 
\begin{equation*}
\beta _{i}=\beta ^{p_{i}},\text{ \ }p_{i}=r_{i0}+r_{i1}x+...+r_{i\left(
j-1\right) }x^{j-1}\in \mathbb{Z}_{m}\left[ x\right] ,\text{ }
\end{equation*}%
and 
\begin{eqnarray*}
\beta ^{m} &=&(\beta _{1}\beta _{2}...\beta _{m})^{(1)}, \\
\beta _{1}\beta _{2}...\beta _{m} &=&\beta ^{p_{1}+...+p_{m}}, \\
p_{1}+...+p_{m} &=&q.x^{j-1}
\end{eqnarray*}%
where $q$ is an invertible element of $\mathbb{Z}_{m}\left[ \left[ x\right] %
\right] $.

\begin{proposition}
The element $\beta =(\beta _{1},\beta _{2},...,\beta _{m})\sigma $ is
conjugate in $\mathcal{A}_{m}$ to $\alpha =(e,...,e,\alpha ^{x^{j-1}})\sigma 
$.

\begin{proof}
Let $h=(h_{1},h_{2},...,h_{m})$ be an automorphism of the tree. Then 
\begin{equation*}
\beta ^{h}=(h_{1}^{-1}\beta _{1}h_{2},h_{2}^{-1}\beta
_{2}h_{3},...,h_{m}^{-1}\beta _{m}h_{1})\sigma \text{.}
\end{equation*}%
Therefore $\beta ^{h}=\alpha $ holds if and only if%
\begin{equation*}
h_{2}=\beta _{1}^{-1}h_{1},\text{ }h_{3}=\beta _{2}^{-1}h_{2},...,\text{ }%
h_{m}=\beta _{m-1}^{-1}h_{m-1},\text{ }h_{1}=\beta _{m}^{-1}h_{m}\alpha
^{x^{j-1}}\text{.}
\end{equation*}%
These conditions can be rewritten as%
\begin{eqnarray*}
h_{2} &=&\beta _{1}^{-1}h_{1},\text{ }h_{3}=\beta _{2}^{-1}\beta
_{1}^{-1}h_{1},...,~h_{m}=\beta _{m-1}^{-1}...\beta _{1}^{-1}h_{1}, \\
h_{1} &=&\beta _{m}^{-1}\beta _{m-1}^{-1}...\beta _{1}^{-1}h_{1}\alpha
^{x^{j-1}}\text{, }
\end{eqnarray*}%
or as%
\begin{eqnarray*}
h &=&(h_{1},\beta _{1}^{-1}h_{1},\beta _{2}^{-1}\beta
_{1}^{-1}h_{1}...,\beta _{m-1}^{-1}...\beta _{1}^{-1}h_{1}) \\
&=&(e,\text{ }\beta _{1}^{-1},\text{ }\beta _{2}^{-1}\beta _{1}^{-1}...,%
\text{ }\beta _{m-1}^{-1}...\beta _{1}^{-1})\left( h_{1}\right) ^{\left(
1\right) }\text{,}
\end{eqnarray*}%
and%
\begin{equation*}
\left( \beta _{1}\beta _{2}...\beta _{m}\right) ^{h_{1}}=\alpha ^{x^{j-1}}%
\text{.}
\end{equation*}%
Since%
\begin{equation*}
\beta _{1}\beta _{2}...\beta _{m}=\beta ^{q.x^{j-1}}\text{,}
\end{equation*}%
we repeat the above procedure replacing $\beta $ by $\beta ^{q}$ and
replacing $h_{1}$ by $\left( h_{1}^{\prime }\right) ^{x^{j-1}}$. This leads
to the conjugation equation 
\begin{equation*}
\left( \beta ^{q}\right) ^{h_{1}^{\prime }}=\alpha \text{.}
\end{equation*}%
In this manner, we determine an automorphism $h$ of the tree which effects
the required conjugation%
\begin{equation*}
\beta ^{h}=\alpha \text{.}
\end{equation*}
\end{proof}
\end{proposition}

\begin{example}
Let $\beta =(e,\beta ^{q})\sigma $ where $q=1+x$. Then $\beta $ is conjugate
to the adding machine $\alpha =\left( e,\alpha \right) \sigma $. Note that
from Example 1, $\beta $ is not obtainable from $\alpha $ by simply choosing
a different transversal. To exhibit the conjugator $h:\beta \rightarrow
\alpha $, constructed in the proof, define the polynomial sequences 
\begin{eqnarray*}
c_{0} &=&1,\text{ \ }c_{1}=q,\text{ \ }c_{n}=2c_{n-2}+c_{n-1}; \\
c_{-1}^{\prime } &=&0,\text{ \ }c_{0}^{\prime }=0,\text{ \ }c_{n}^{\prime
}=c_{n-1}+c_{n-1}^{\prime }\text{.}
\end{eqnarray*}%
Then
\end{example}

\begin{equation*}
h=\left( e,e\right) ^{\left( 0\right) }\left( e,\beta ^{-1}\right) ^{\left(
1\right) }\left( e,\beta ^{-\left( 1+q\right) }\right) ^{\left( 2\right)
}...\left( e,\beta ^{-c_{n}^{\prime }}\right) ^{\left( n\right) }...\text{.}
\end{equation*}

\textbf{Step 4.} By Proposition 2, we have $A\leq \overline{A}=C_{\mathcal{A}%
}\left( \alpha \right) $ and%
\begin{equation*}
A^{h}\leq C_{\mathcal{A}}\left( \alpha ^{h}\right) =C_{\mathcal{A}}\left(
\beta \right) =\overline{\mathit{D}_{m}\left( j\right) }\text{.}
\end{equation*}%
This finishes the proof of the theorem.

\end{document}